# A Conversation with Chris Heyde

**Paul Glasserman and Steven Kou**

*Abstract.* Born in Sydney, Australia, on April 20, 1939, Chris Heyde shifted his interest from sport to mathematics thanks to inspiration from a schoolteacher. After earning an M.Sc. degree from the University of Sydney and a Ph.D. from the Australian National University (ANU), he began his academic career in the United States at Michigan State University, and then in the United Kingdom at the University of Sheffield and the University of Manchester. In 1968, Chris moved back to Australia to teach at ANU until 1975, when he joined CSIRO, where he was Acting Chief of the Division of Mathematics and Statistics. From 1983 to 1986, he was a Professor and Chairman of the Department of Statistics at the University of Melbourne. Chris then returned to ANU to become the Head of the Statistics Department, and later the Foundation Dean of the School of Mathematical Sciences (now the Mathematical Sciences Institute). Since 1993, he has also spent one semester each year teaching at the Department of Statistics, Columbia University, and has been the director of the Center for Applied Probability at Columbia University since its creation in 1993.

Chris has been honored worldwide for his contributions in probability, statistics and the history of statistics. He is a Fellow of the International Statistical Institute and the Institute of Mathematical Statistics, and he is one of three people to be a member of both the Australian Academy of Science and the Australian Academy of Social Sciences. In 2003, he received the Order of Australia from the Australian government. He has been awarded the Pitman Medal and the Hannan Medal. Chris was conferred a D.Sc. *honoris causa* by University of Sydney in 1998.

Chris has been very active in serving the statistical community, including as the Vice President of the International Statistical Institute, President of the Bernoulli Society and Vice President of the Australian Mathematical Society. He has served on numerous editorial boards, most notably as Editor of *Stochastic Processes and Their Applications* from 1983 to 1989, and as Editor-in-Chief of *Journal of Applied Probability* and *Advances in Applied Probability* since 1990.

His research has spanned almost all areas of probability and statistics, ranging from random walks to branching processes, from martingales to quasi-likelihood inference, from genetics to option pricing, from queueing theory to long-range dependence. He has edited twelve books, and authored or co-authored three books, *I. J. Bienaymé: Statistical Theory Anticipated* (1977), with E. Seneta, *Martingale Limit Theory and Its Application* (1980), with P. Hall, and *Quasi-Likelihood and Its Applications* (1977). Chris Heyde has been an outstanding citizen and leader of the probability and statistics research community.

The interview ranges over his education in Australia, moves to the USA and the UK, return to Australia, appointment at Columbia, major research contributions, and professional society and editorial activities. It ends with a look forward in time and some concerned comments about the future for statistics departments.

*Paul Glasserman is Jack R. Anderson Professor of Business and Senior Vice Dean, Graduate School of Business, 101 Uris Hall, Columbia University, New York, New York 10027, USA e-mail: pg20@columbia.edu. Steven Kou is Associate Professor, Department of Industrial Engineering and Operations Research, 313 Mudd Building, Columbia University,*

*New York, New York 10027, USA e-mail: sk75@columbia.edu.*

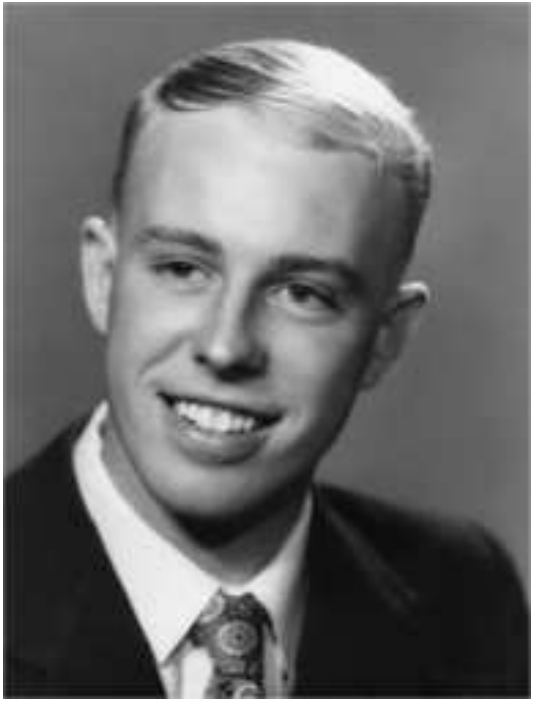

Fig. 1. *Chris Heyde as a student at the University of Sydney (1957).*

The interview took place on November 24, 2003, at Chris Heyde's office.

**Paul:** Chris, tell us about growing up and studying in Australia. Why did you do your doctoral work in Australia when most similar students did their graduate studying overseas?

**Chris:** Australia has been described as the lucky country. It certainly was thought of in those terms 40 years ago (e.g., [8]) and largely still is. It's very untroubled and offers a high quality lifestyle. And I feel fortunate to have grown up there and to have had much of my life there. I had a schooling that was very unpressured and I was mostly interested in sports, at which I was pretty good. I never thought much about things that were academic, until the age of 16 when I pulled my Achilles tendon, and I was unable to continue with athletics. So I got interested in my studies for the first time. I discovered that the mathematics teacher was rather stimulating, and it wasn't too late to catch up at that stage. I had two years of high school still to go and managed to finish up as dux (top student) of the school. That embarked me on an academic career.

The natural thing for someone living in Sydney was to go to the University of Sydney as an undergraduate. And the education system there was such that if one did a science degree, and that was the automatic thing for me, you started off with four first-year subjects and followed these with three second-year subjects, two third-year subjects and one fourth-year subject. So you could follow a path of least resistance through the sciences, dropping the least favored subject at the end of each year. I started off the first year with mathematics, physics, geology and chemistry and I could have ended up majoring in any of them. But mathematics first and foremost caught my fancy. In my fourth year, it was the first year of the appointment of a new Professor of Mathematical Statistics—Oliver Lancaster. He had trained as a medico and had worked in medical as well as mathematical statistics. It was the first time there was an option to do Honors in Mathematical Statistics rather than Pure Mathematics or Applied Mathematics. And I had some difficulty in deciding what I wanted to do as I had nearly decided to go into seismology, but eventually Mathematical Statistics won and I am glad in retrospect that I took that decision. So I did Honors in Mathematical Statistics and then it seemed natural to stay on and do a Master's degree. Even if I had decided to go overseas, the overseas academic year starts in September whereas the Australian academic year starts at the beginning of the year so you can virtually finish a Master's degree before going overseas. But while I was doing a Master's degree Oliver Lancaster took me to Canberra—to the Australian National University (ANU)—where he was having discussions with Pat Moran, who was the Professor of Statistics there in the Institute of Advanced Studies. I was attracted by the environment at ANU and the fact that the Institute at ANU was focused on research and there were no undergraduates; there were only Ph.D. students and faculty. It was quite a remarkable scholarly environment and I could see that exciting things were happening. At that time applied probability was growing into an identifiable discipline and it was visibly happening in Canberra. And ANU seemed to be a nice place to live as well, with a college environment (University House) where all the single Ph.D. students (of whom some 50% were from outside Australia) were living, as well as quite a few of the faculty, and this was just a few minutes walk away from work. So, I thought it was a good idea to continue there rather than to go off to Cambridge as had usually been the case for the best people graduating out of the University of Sydney. I think, in retrospect, that it was a good decision. Lots of interesting science was going on and living in University House was pleasant socially and culturally. Indeed, I met my wife-to-be in University House. So I stayed there to finish my Ph.D.

**Steve:** What did you study at the Australian National University?

**Chris:** Well, when I was doing my Master's degree at the University of Sydney, I was working on moment problems—in particular, when is a distribution determined by its moments? It wasn't that anyone locally had any particular expertise in the area, but I just happened to get interested in it and so I did it for a Master's degree. A number of papers came out of that research and the one that is most remembered is the one that shows that the log-normal distribution is not determined by its moments [4]. Then when I went to Canberra to start my Ph.D. Pat Moran said to me, in a very brief conversation, "there are lots of interesting issues associated with what you can tell about a distribution on the basis of its passage time properties. Why don't you have a look at that?" So essentially I went away for three years and wrote a thesis on that. No one took any real notice of what I was doing. The strategy at that particular stage was of benevolent neglect of the students who were not working on collaborative projects. The faculty would talk to the students if it turned out to happen, but otherwise they just left them to their own devices and they either survived or didn't survive on their capacity to work independently. But it all worked out alright for me and quite a few papers came out of my Ph.D. research—ultimately enough to get me a job overseas which was the natural thing to do after these years in Australia. It was recognized if you didn't do a Ph.D. overseas then you would have to go and work overseas even if you were going to come back to Australia—the international experience was essential.

**Paul:** When did you decide that you wanted to pursue an academic career?

**Chris:** It almost happened automatically and inevitably at the stage when I did a Ph.D. and there were relatively few industry jobs for people with this sort of training, and it was clear at that stage that there were going to be plenty of academic jobs—that was the beginning of the 1960s. The 1960s were a very expansionary period in academia and there were new universities being created all over the place and there were so many positions about, that if in fact I had stayed in Australia, I probably could have been appointed immediately as a senior lecturer upon getting my Ph.D. There were just so many vacancies but I knew I should go overseas and get the experience overseas at that stage.

## TO THE USA AND THE UK

**Steve:** What was your first position after completing your Ph.D.?

**Chris:** I went to Michigan State University (MSU) and I did so because Joe Gani, who had been a faculty member at ANU in the Institute of Advanced Studies while I was studying, had taken a position there as had Uma Prabhu from the University of Western Australia whom I also knew quite well. Uma Prabhu had invited me to give a lecture course at the University of Western Australia while I was still a Ph.D. student at ANU. We were looking at the possibility of creating a new stream in stochastic processes at MSU. Joe Gani had recently started *Journal of Applied Probability* in Australia and this was the beginning of an autonomous literature in applied probability. It was clear that the time had come for this—the subject had grown to a stage where it was a recognizable discipline in its own right. Anyway, I finished my Ph.D. in August 1964 and then went to East Lansing.

**Steve:** How long did you stay in East Lansing?

**Chris:** Just a year. It became quite clear to us that we were not easily going to be able to start the new program in East Lansing. Joe Gani decided to take the Chair of the Statistics Department at the University of Sheffield in Britain, Uma Prabhu went to Cornell and I decided to go with Joe Gani to Sheffield.

**Paul:** How long were you in Sheffield?

**Chris:** Two years in Sheffield and one year in Manchester. During the time that I was in Sheffield, Manchester statistics, which had had a distinguished history first under Maurice Bartlett and then under Peter Whittle, fell apart with most of its faculty leaving to go to other institutions. Joe Gani had the opportunity to create a Manchester-Sheffield School of Probability and Statistics. And I was sent across to take charge of the Manchester component, which I did for a year. But during that year I had three offers of jobs in Australia: two of them professorships—head-of-department type jobs, one of them a readership back at ANU and I decided to take the readership. I thought that I should take one of these jobs because academic jobs at the senior level in Australia had not tended to come up frequently, and once a professorship was filled it was often occupied for 20 or 25 years. Those were the days when there was only one (full) professor in a department and it was impossible to be promoted internally through to (full) professor. You had to move in order to get a professorship and there were only a small number of these. Also at this stage my wife and I had our first child and we were attracted at the idea of bringing up our children in Australia. So that was why we decided to return at that stage.

**Paul:** That was your son…?

**Chris:** I have two sons, yes.

**Paul:** Named?

**Chris:** Neil and Eric. Neil was one year old when we left England. He is now resident in London and is a musician and a faculty member at the Royal Academy of Music. My younger son Eric has a Ph.D. in Electrical Engineering and he lives in Sydney.

**Paul:** I wanted to ask about your research interests during your time in Sheffield and Manchester.

**Chris:** By the time I finished my Ph.D. I had developed interests in properties of random walks, in particular the characteristics of sums of independent and identically distributed random variables as the sample size grows. So initially I continued exploring these things. I was very much concerned with associated questions like what happens asymptotically to renewal functions. Anything to do with properties of sums of independent random variables took my fancy and I got interested in particular in rate of convergence results associated with the central limit theorem, the laws of large numbers and the iterated logarithm law. So I was exploring all these things and it was not 'til I went back to Australia that I started to get involved with a broader range of stochastic models with dependent random variables.

## RETURN TO AUSTRALIA

**Steve:** It seems that you've returned to Australia many times.

**Chris:** Well yes. First I went back to ANU in Canberra. ANU is an unusual institution in being the amalgam of two separate institutions. The Institute of Advanced Studies, where I had done my Ph.D., was created immediately post-World War II, as a research-only institution, essentially with the purpose of training the next generation of researchers and academics for the Australian university system. In Canberra, however, there had also been what was originally called the Canberra University College, a teaching institution that was created in the 1930s and was originally under the administration of the University of Melbourne. In about 1960 the then Prime Minister, Robert Menzies, "married" these two institutions more or less over the dead bodies of both of them and the marriage really was never consummated. So there were these two institutions, the Institute (the original ANU) and the Faculties (which had been the Canberra University College) on the campus, largely separated physically and philosophically. Fortunately the Statistics departments spoke to one another but many of the other departments of the university were not in good communication. The Statistics departments interacted well because Ted Hannan, who was then the Professor in the Faculties department, had been one of the two first Ph.D. students of Pat Moran—the other having been Joe Gani. Anyway, I went to the Faculties department at the end of 1968 and stayed there 'til 1975.

In 1974 Joe Gani had been asked to come out to Australia by CSIRO—Commonwealth Scientific and Industrial Research Organization—which is the government scientific organization. CSIRO is a countrywide organization which had approximately 50 scientific divisions, including one then called Mathematical Statistics which had some 50 staff. The Chief of this division had recently died, and CSIRO was looking to reorganize and perhaps to expand the division. They asked Joe Gani to come from the UK and advise them on the possibilities and he wrote very detailed recommendations. They then accepted these and asked him to come and implement the program. So although he had been happy with the Manchester–Sheffield School he decided to take the job, in part because he wanted to repay the country which had accepted him as an immigrant in 1948. He also very much wanted me to join him. He had a mandate to substantially build up the division which had principally been a consulting group providing statistical advice within the organization, and doing a certain amount of basic research, but not basic research first and foremost. The new mandate was for a much broader range of activities and to put the division on the world stage. So I found this attractive and went to join him in 1975 and I was there ultimately 'til the beginning of 1983. It was a period of considerable success.

Joe Gani was appointed for a seven-year term and towards the end of that period politics had intruded into the organization and there were changes of policy in the wind. The organization had started to feel it necessary to become more inward looking and devote much more of its time and energy to doing science aimed at the economic benefit of Australia. So Joe Gani didn't want to stay and run a division with much reduced scope. He went off to the United States and I stayed as Acting Chief of the division for a period 'til I could get a good university chair.

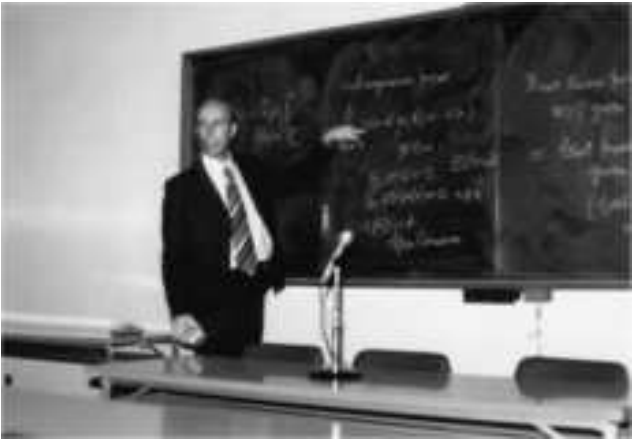

FIG. 2. *Nagoya University, Japan, 1977. Seminar by Chris Heyde.*

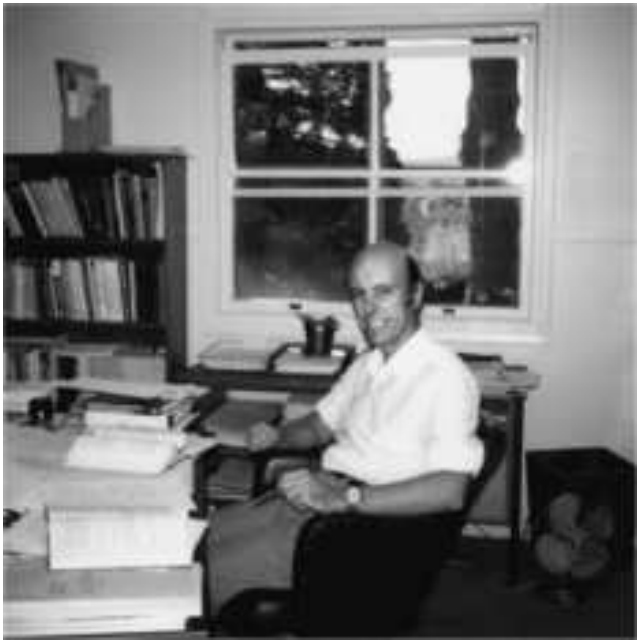

FIG. 3. *Chris Heyde in his office at CSIRO, Canberra, 1979–1980.*

Fortunately the position of Head of the Department of Statistics at the University of Melbourne came up and I went there. I very much enjoyed my time in Melbourne. The then Vice Chancellor was very supportive of what I wanted to do, such as to start a statistical consulting center run as a commercial operation. The university was willing to pay for the salary of a director for the consulting center and to refurbish a set of offices to commercial standard. The consulting center had to pay for everything else itself out of its generated revenue. But it was perfectly able to do that, and it has subsequently been a real success. The Vice Chancellor also gave me some extra positions so I was able to build quite a lively environment. In addition, this was the time when the Australian Government was providing money for what were called Key Centers and I was fortunate in being able to collaborate successfully with three other universities in the Melbourne metropolitan area to form a Key Center for Statistical Science. It was the University of Melbourne together with Monash University, La Trobe University and the Royal Melbourne Institute of Technology that banded together to form this center, of which I was the Foundation Director. Collectively we were a very strong group. Individually we had small numbers of honors students and masters students—in fact too small to run a broad spectrum of courses—but collectively we could do it. The students used to travel from one institution to another depending on the day of the week and they had access to a very broad range of topics. This was very healthy and it's still an operating arrangement more than 20 years later. It is necessary in a place like Australia, where the individual universities have relatively small numbers of students, to develop cooperative schemes. Of course any such schemes have difficulties and there always are rivalries between institutions. In setting up the Key Center it was somewhat amusing, albeit frustrating, that the principal obstacles were in getting each institution to accredit the courses from the others for the purpose of their degree. Academics, of course, like to look at other people's courses and criticize them, and to get agreement on this sort of thing did take some time and effort. But ultimately we did it.

So I was at the University of Melbourne, things were going well, and I thought I would probably stay there for my whole career. But then the position of Professor of Statistics in the Institute of Advanced Studies at ANU came up—and I was successful in getting it. So I ended up going back to where I had been a Ph.D. student, and I've essentially been there ever since, but of course I have had a Columbia connection since 1993.

## THE COLUMBIA CONNECTION

**Paul:** How did the arrangement with Columbia come about?

**Chris:** Well, in 1990 I had become Editor-in-Chief of *Journal of Applied Probability* and *Advances in Applied Probability*, and it was my policy to take the opportunity to visit people who were associate editors of the journals if it was convenient in my travels. I came to Columbia for a few days on such a visit. Now the Columbia Statistics Department had been in some difficulty after key retirements and resignations and had nearly been closed down. They were anxiously hoping to rebuild, and they strongly suggested that I come for a semester as a visitor. Although I was a bit frightened by New York in those days, I agreed. So I came in 1992 to spend the fall semester at Columbia and I very much enjoyed the experience. I didn't find New York as I had feared and I had some very good students in my courses—one of whom is taking part in this interview. So I was pleased with the experience. Then Columbia essentially made me an offer that was too good to refuse at the time, and I've been coming back ever since. Just putting this in context, if you go back to 1992, I was just coming to the end of a three-year term as Dean of the School of Mathematical Sciences at ANU and that had been a difficult assignment. The School of Mathematical Sciences was the first group in the University to bridge the Institute/Faculties divide. I did mention earlier that marriage between the Institute and the Faculties had never been consummated. Well, it was the School of Mathematical Sciences that first joined Institute and Faculties components of the University. There were plenty of associated troubles and the thought of the comparative freedom of Columbia was particularly attractive. Also I had the sense of things being achievable in the United States that were not achievable in Australia at the time. In a sense the pioneering spirit is still alive and well in the United States and Columbia at the time was very fortunate to have a Vice Provost, Michael Crow, who was very entrepreneurial. Whenever one went to him to discuss things that might be done the response was always to seek a way of doing it. This was such a contrast to what I had been finding in recent times in Australia where similar visits to administrators met with basic negativity. It was very attractive to me to have the opportunity to do things—for example, get the Center for Applied Probability launched at Columbia. There was strong support from the University helping to make us competitive so that we were ultimately able to get a General Infrastructure Grant from the National Science Foundation. So this was a breath of fresh air for me after the comparative difficulties of the previous three years.

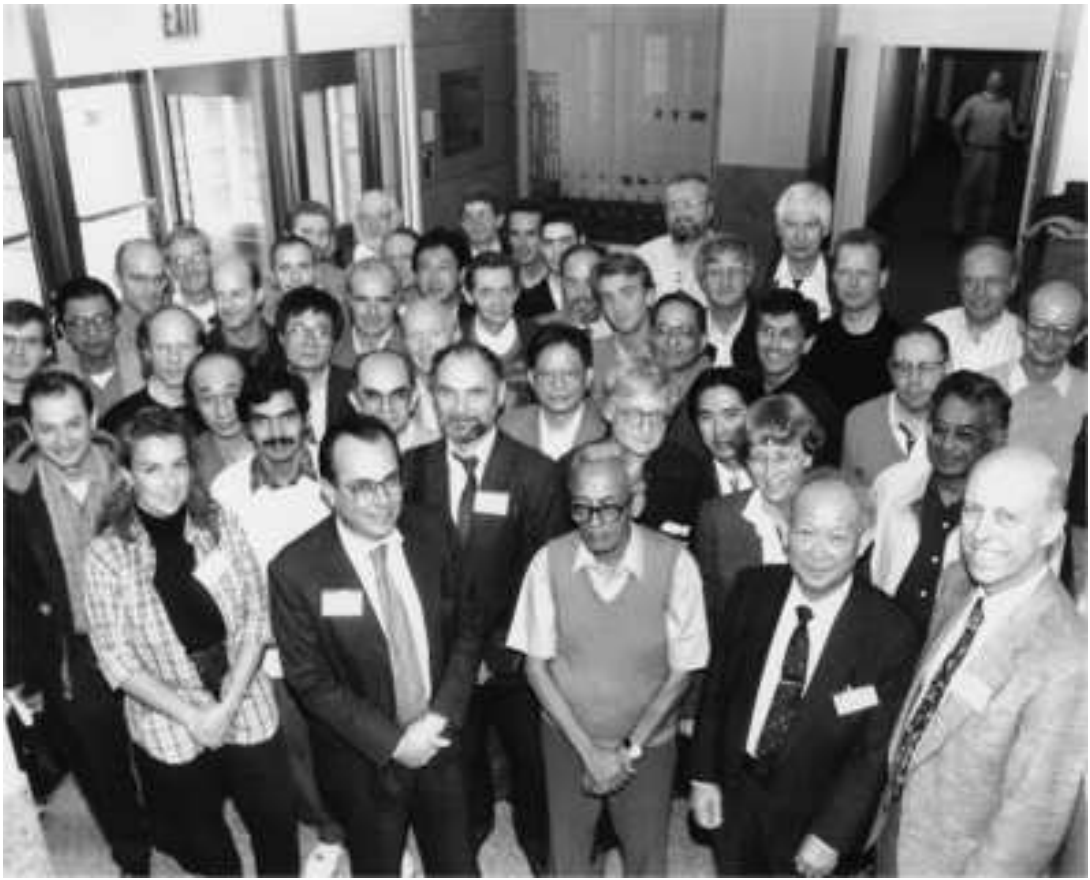

Fig. 4. *Probability Towards 2000 Conference, Columbia University, New York, October 1995. A distinguished collection. Organizers L. Accardi (front row, second from left) and C. Heyde.*

## SOME RESEARCH CONTRIBUTIONS

**Steve:** Early in your career your research took an important turn from a focus on independent random variables to models with dependence.

**Chris:** Yes. While I was a Ph.D. student the faculty who were at ANU were much interested in Markov processes and so of course I absorbed the culture of Markov processes. But I had always felt that this was just one step away from independence and that if one really wanted to have a broadly applicable tool one had to develop something rather more general. But from the early 1960s, which was of course the time I was doing my Ph.D., there were beginning to be a few papers on sums of dependent variables. I was ultimately quite influenced by the book that Patrick Billingsly wrote in 1961 on inference for Markov processes. Although he was dealing with stationary ergodic Markov processes, he used the martingale central limit theorem as a basis for doing parameter estimation. I saw this and realized that martingales were potentially a very general tool. You can actually make martingales out of any integrable stochastic process. If you take off the conditional means and you sum up, you have a martingale. So I realized that if there were good limit theorems for martingales, it would be a vehicle for doing a whole lot of investigative modeling work—and most of those good limit theorems were able to be developed in the 1970s.

An interest in all sorts of stochastic models had originated from my time as a Ph.D. student in Canberra. I realized that if you're going to study realistic and useful probability models you would also need to be able to do the statistical inference for the models.

The modeling philosophy at ANU probably came about because Pat Moran had been an Australian Scientific Liaison Officer in London during World War II. This had required him to travel around military research establishments seeing what research was being done and reporting on it to Australia. And the experience had given him a very broad overview of quantitative methods in science which he took into his own research. Although I didn't work directly with any of the faculty, the culture of mathematical modeling practiced in the department rubbed off on me, and it strongly influenced the way I've subsequently thought about science. So it was fairly natural that I would start looking at a wide variety of different stochastic models and studying them and their associated model validation and inference.

**Paul:** Your book on martingales with Peter Hall [1] has had a major and lasting impact. And yet some of these tools are not as well known as they could be.

**Chris:** I think it's very regrettable that statisticians have not adopted martingales as a tool nearly as widely as might have been expected. It's largely in the United States rather than in other parts of the world where this has happened. It seems to come about because in the United States there is a notional division between probability and statistics, so that students who specialize in statistics will typically not get schooled in things such as martingales or dependent variables. Only students who specialize in probability will tend to get such courses. The consequence is that statisticians with the standard training mostly don't know about martingales and they don't realize that a lot of the inferential things that they do can be usefully generalized without any real additional cost. This is a particular characteristic of the United States system. In many European countries and certainly in Australia there's no divide between probability and statistics; it's all regarded as part of one continuous spectrum of statistical science activity.

**Steve:** Your book on quasi-likelihood [6] seems to be a bridge between probability and statistics. Do you agree?

**Chris:** Yes, this was a natural development from the earlier work. I had written about likelihood-based methods in the martingales book and explored general properties of the maximum likelihood estimator. Martingales are certainly the natural way of looking at the maximum likelihood estimator because the derivative of the log-likelihood is a martingale almost universally. So martingales then provide the basis for studying the asymptotics from which all the large sample properties of the estimator come. It was fairly natural to try to do inference in a broader setting. I had always been quite concerned about the strength of the assumptions that had to be made in order to use the full maximum likelihood theory. So it was attractive to investigate results that conveyed much of the advantage of the likelihood theory, but didn't require the full knowledge of the distribution, just assumptions about first and second moments—a covariance structure. From the mid-1980s I had been thinking about these issues and talking about them with people like V. P. Godambe from Waterloo. Ultimately this led on to the book, and although my motivation was general theory of inference I had stochastic processes in the back of my mind. The theory is applicable to sample surveys and random fields and indeed any sort of stochastic system.

**Paul:** What do you view as some of the other landmarks in your research career?

**Chris:** I guess it was rates of convergence work that established my reputation. Then came the martingales work followed by the statistical inference. I've also done some significant things for various stochastic models, branching processes being one. For example, there are particular constants called the Seneta–Heyde constants that can be used to normalize the supercritical Bienaymé–Galton–Watson process to obtain almost sure convergence. I have also done some important work on time series and made useful contributions to the history of statistics.

**Paul:** The time series work was with Ted Hannan?

**Chris:** Yes, that goes right back to the start of the martingales story. Ted Hannan was a time series expert and he published a very influential book in 1970 [2]. While he was writing this he used to come to my office every day to talk about time series and at this stage I was beginning to write a review paper about martingales which was published in 1972 [5]. Ted and I started to think about what the martingale property means for time series innovations and we wrote a paper in 1972 [3] that showed that the best linear predictor is the best predictor if and only if the innovations are martingale differences. This is a very important result because it says that you would use a linear model if the innovations are martingale differences, but if not it would be wise to look for a better nonlinear model. If you can find one you may have much better prediction error properties. So this was important for time series, which had very much grown up in the mold of the linear theory because it came from Gaussian process theory where you can quite satisfactorily represent everything as a linear model with independent innovations because of the Wold decomposition theorem and associated results.

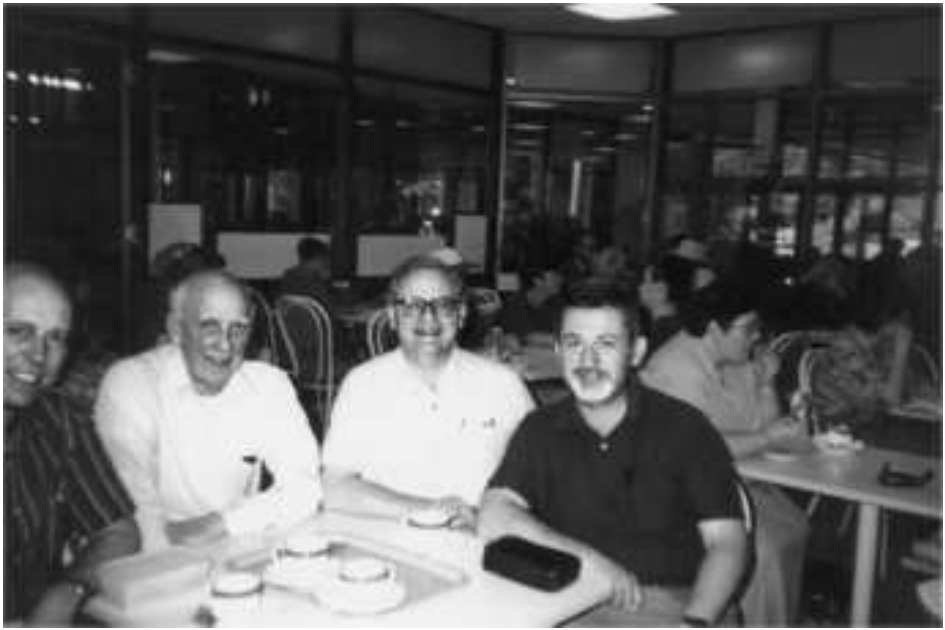

FIG. 5. *Canberra, 1994. Chris Heyde, Ted Hannan, Joe Gani, Eugene Seneta.*

**Steve:** Tell us a little more about your research on branching processes and in the history of statistics.

**Chris:** My interests in branching processes were kindled in the late 1960s by my colleague Eugene Seneta, who had started working in the area some years earlier. I was attracted by the beautiful and then nonstandard limit theorems that could be obtained, many involving sums of random numbers of random variables and also random normings. Through working in the area I ultimately became involved with a diverse collection of problems in biology.

It was also with Eugene Seneta, and through our then common interest in branching processes, that we became involved in the history of probability and statistics. This began with investigating the origins of the study of population models and expanded into an examination of nineteenth century French probability and statistics. History has been a long-term interest for me, but rather more in the nature of a hobby than of a primary research theme.

**Steve:** You have managed to work comfortably in both probability and statistics. Do you think that's something more researchers should try to do?

**Chris:** It seems to me that nowadays there is a tendency towards overspecialization. I think it's very important that students take a broad view of their subject and it's crucial that university faculty promote such a broad view. We are training students for careers that will go over 30 or 40 years during which time it will be necessary for them to react to diversely changing circumstances. So we should equip them as best we can to develop new methodology that goes along with the changing circumstances. It's certain that the focus of applications will change over time and that people will need to be quite flexible in their capacity to approach problems. This will be the case whatever one's area of specialization. So I think it's most unwise not to give a broad-spectrum core curriculum to students, to allow them a good springboard for change.

**Paul:** Do you consider statistical consulting an important part of graduate education in statistics and do you consider it an important activity for statistics departments?

**Chris:** I do. I had quite an involvement with statistical consulting, although principally at the management level.

In CSIRO I used to direct the activities of quite a large group of people who were mostly doing consulting. We used to sit around the table at least once a week and discuss the problems that came in and the work that was going to be done. Also while I was in CSIRO we started, out of the division, a commercial consulting company called SIROMATH. CSIRO was one of three shareholders, the other two being commercial organizations. Then when I went to Melbourne I started the University of Melbourne Statistical Consulting Centre. I had always very much seen value in the stimulus provided to statistics departments by the ongoing consulting. It really changes the sort of discussions that one has around morning or afternoon tea/coffee time. If there was no statistical consulting going on in the department the conversation might be about sport. But if there was consulting going on in the department, it was usually about consulting, and I found that to be a considerable stimulus to both the staff and the students. I think all students need exposure to these experiences and some sort of apprenticeship to help them get a start. Otherwise they find that bridging the gap between the theory they've learned in class and what they actually have to do in real-world applications is very hard.

## PROFESSIONAL CONTRIBUTIONS

**Paul:** You've long been a member of the Australian Academy of Science and you were just recently named to the Australian Academy of Social Science as well; now what can you tell us about that?

**Chris:** Well, I became a Fellow of the Australian Academy of Science in 1977 and I served the academy in a number of roles over the years. I was a Council member, and the Academy Treasurer, and the Academy Vice President, and chair of the committee that set up the Australian Foundation for Science. I've always been interested in the broader issues of science in society. Recently I became a Fellow of the Academy of Social Sciences as you've mentioned, and I'm one of only three people who's a Fellow of both those academies. So I hope that I can do something to foster collaboration between these academies. There is a tendency for such organizations to go their own way and they can potentially achieve more if they work cooperatively—especially since the governments don't always take much notice of the propositions that are put to them. The larger the constituency you have, the more likely is your success.

I guess I qualified for the Academy of Social Sciences because my work on martingales and time series has had a significant impact in econometrics, and also because of my work on the history of statistics. Curiously, martingales seem to have created more interest in econometrics than in statistics.

As far as history is concerned, I think it enriches the study of every subject and I try to give a historical perspective to students. If you look for the motivation behind most of the really important discoveries, you typically find it to be a very practical problem of the time. This should enhance our appreciation of the exposure to good applications that we are so fortunate to have in the statistical sciences.

**Steve:** Throughout your career you've been active in advancing the research community through centers, professional societies and editorial work. This must have taken time away from your research, so you must feel that these activities are important.

**Chris:** Well, I've taken the view that contributions to the well-being of one's profession are an integral part of one's professional life. I think that our profession is not well understood or recognized by the community and that to a significant extent it's our fault for the failure to successfully promote ourselves. Statistics does have a name with unfortunate connotations, because many people see the subject as intrinsically dull and boring, and even misrepresentational. So we have to try to confront those negatives and replace them with a more positive image. I have been concerned with promoting the professional societies in statistics, trying to establish a voice for the statistical profession and, indeed, a seat on the relevant committee whenever there has been an official inquiry where a statistical input was important.

One such experience I had was when I was a member of the Australian Government's Scientific Advisory Committee to review the possible effects of herbicides and pesticides on veterans who served in Vietnam during the Vietnam War. There had been quite a lot of complaints from veterans about disabilities that they and their families seemed to be suffering. Both the Australian and the United States governments were essentially obliged to have large-scale investigations. Ultimately there were something in excess of 80 different disabilities that were investigated, only one of which (chloracne) could be categorically attributed to exposure to the herbicides. The others were more general disabilities and, rather strikingly, no significantly statistical results appeared in any of this. There were many complex statistical issues that occurred in this investigation. I think it is very important that the profession seeks these tasks and that not doing so contributes to our lack of identity.

**Paul:** What are some of the other things you've done through professional societies?

**Chris:** I've always been very concerned to help people realize the value of promoting the societies, nationally and internationally. So I've tried to build up membership, to establish additional branches of societies and to establish a culture of professional participation. One of the things I have been concerned about is a fragmentation of the discipline as one sees many societies creating more and more special-interest sections. It is often the case that people who belong to one special-interest section don't talk to people who belong to other special-interest sections. And so the force which unifies people within the profession is relatively weak and the professional coherence is being dissipated over time. Ultimately this is a considerable danger to the profession.

**Paul:** Tell us about the Bernoulli Society.

**Chris:** The Bernoulli Society got its start as a section of the International Statistical Institute (ISI) in a rather indirect way. There had been a number of international special-interest groups that started up that were essentially autonomous. One of them was a group called the Committee for Conferences on Stochastic Processes which was essentially started by Uma Prabhu, Richard Syski, Julian Kielson and Wim Cohen along with the journal *Stochastic Processes and Their Applications*. They had a series of conferences, which were held every two years and moved around the world. Now there was a certain amount of effort into looking into what might provide an international umbrella for such organizations. This ultimately led to the ISI providing an umbrella through the Bernoulli Society, which was created in 1975. It is the largest section of the ISI and is intended to embrace mathematical statistics and probability.

When the Bernoulli Society was started it needed a mission and it was decided that it should have large international conferences. The first of these was held in Tashkent in 1986. I was president at that time. The conference has subsequently gone on to be held every four years and to move widely around the world. It has been very successful and it is now generally held jointly with the Institute of Mathematical Statistics. Fortunately the international organizations tend to work together more nowadays than used to be the case.

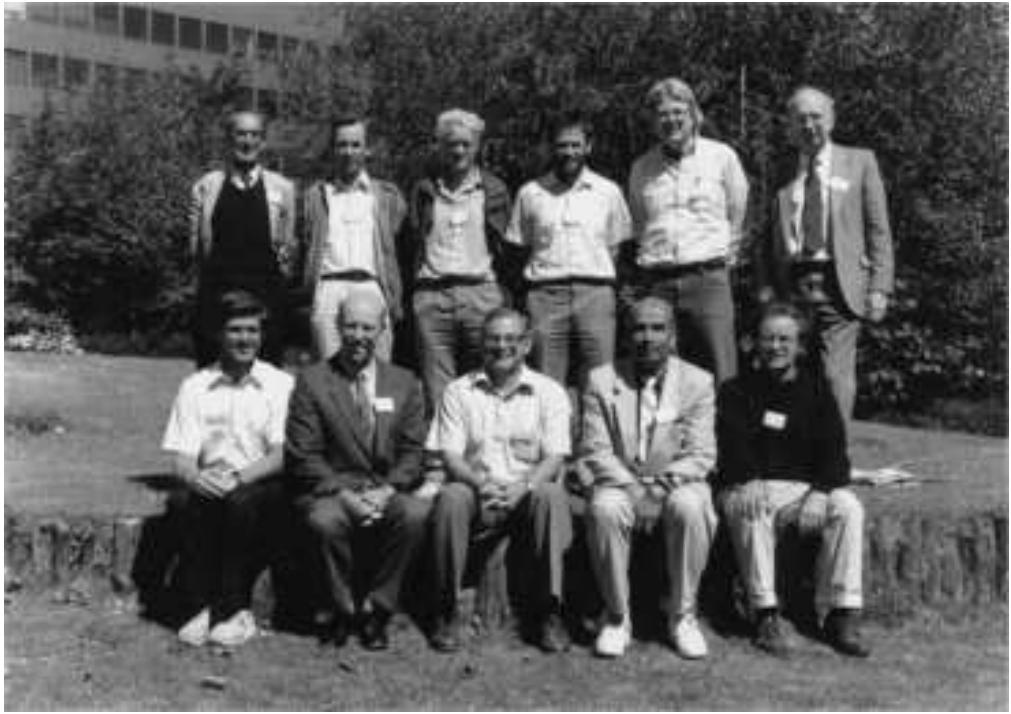

Fig. 6. *Sheffield Conference on Applied Probability, August 1989. Plenary speakers and organizing committee. Front row l. to r.: R. L. Taylor, C. C. Heyde, J. Gani, N. U. Prabhu, C. Anderson. Second row l. to r.: D. G. Kendall, D. Grey, C. Cannings, J. Biggins, R. Durrett, R. M. Loynes.*

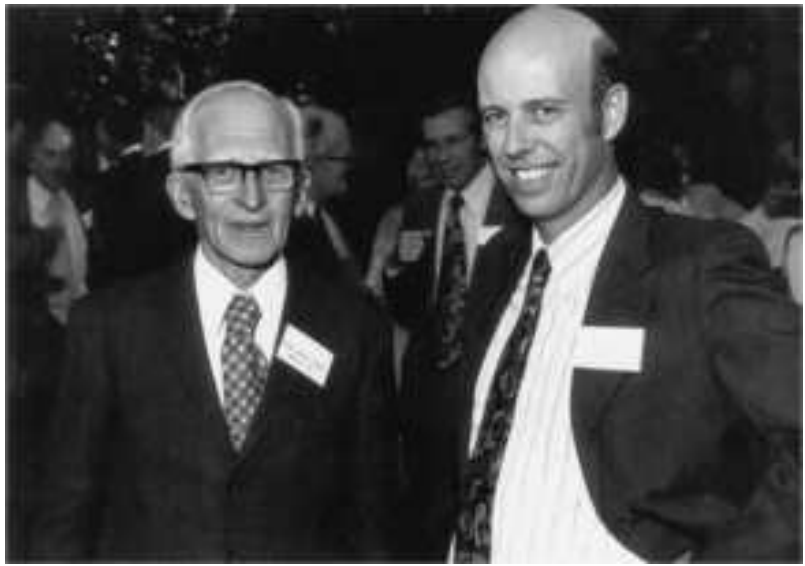

Fig. 7. *Canberra, 1978. Edwin Pitman and Chris Heyde on the occasion of the award of the inaugural Pitman Medal.*

Society membership is generally in decline and every effort needs to be made to arrest the trend. The professional organizations must offer services that are perceived to be of value by their potential members.

**Steve:** During your career you have been much involved with editorial work. What have you learned from this?

**Chris:** I think that perhaps the strongest single message is that people don't write with their audience in mind as often as one might hope. We typically don't teach our students how to write research papers. We take it for granted that they'll know how to do this and that's not really the case at all. People need to market their work. They have to say why it's important and how it relates to what has gone before in as compelling a way as they reasonably can. The idea that science is going to be valued objectively on the basis of its quality without any particular effort being made to promote that quality is quite false. People have to work hard to display their wares convincingly. I have spent a huge amount of time over the years in trying to help people improve their papers to the extent that the audience is at least likely to see what the basic message is. Otherwise there is little point in publishing the work, even if there are gems within it. I think that one of the key roles of an editor is to recognize quality and to avoid the mistake of rejecting obscure papers of real intrinsic value. And particularly with inexperienced authors, you should be sympathetic to their circumstances and try to help them. I've always tried to look at the things from a human point of view as well as an academic point of view but sometimes this involves telling people bluntly that they should have done a better job.

## A LOOK FORWARD

**Paul:** What are your thoughts about the future of probability and statistics?

**Chris:** I think that probability and statistics have a good future in prospect and there's much that's interesting going on. Where probability will have its home in the long run, and under whose auspices, is an open question. In the United States probability is fairly commonly located in mathematical departments, and whether that will disadvantage it in the long run from being as actively involved with applications as might be the case—I'm not sure—but it's potentially a danger. I'm more concerned about the home for statistics. I see statistics having most of its vitality and life in its burgeoning applications and it is clear that the effort in contributing to what you might describe as core methodology is significantly decreasing. This core methodology is that which is applicable to a broad spectrum of statistical applications, and if the subject is diverging into fields which are not in communication, then the context-free methodologies that might be useful in each may not be developed and widely advertised. Central core methodology and its maintenance are very important. Actually, I doubt that there is a long-term future for the statistics department as we know it. I rather suspect that in the long run there will be some other common institutional structure and perhaps a name which indicates a broader spectrum of quantitative endeavor. And an amalgam of special skills over a broader range of the discipline than we would see at present. I think our present arrangements are too limiting and don't take account of the fact that there's far more statistics being practiced outside of statistics departments than inside. Of course it is the same for probability. Risk and chance are everywhere.

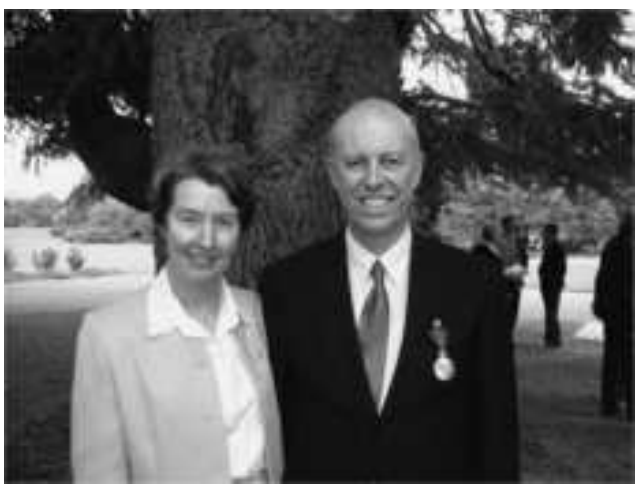

FIG. 8. *Government House, Canberra, April 2003. Chris and Beth Heyde after the award of the Order of Australia.*

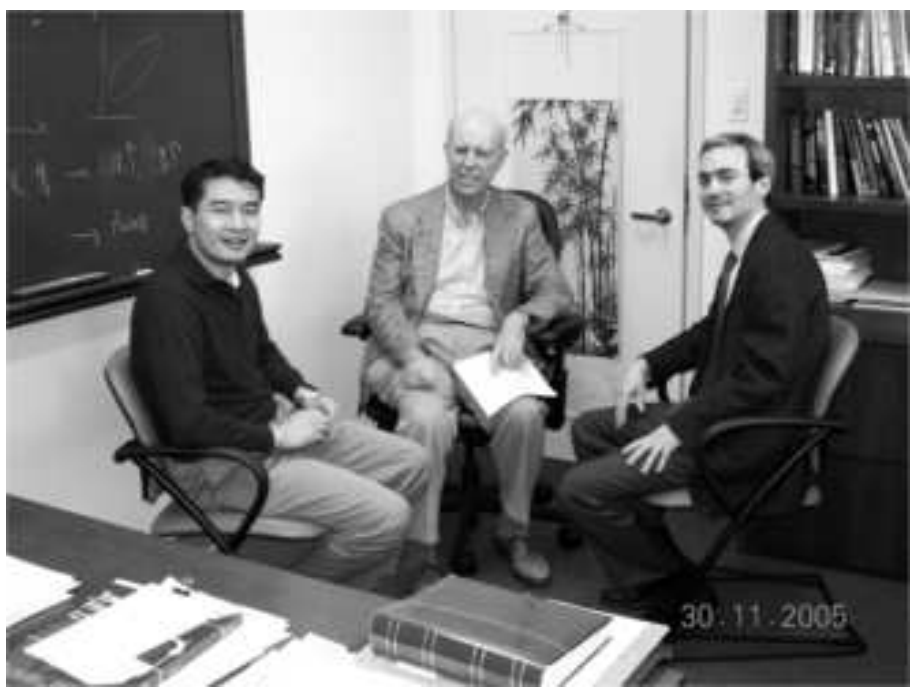

FIG. 9. *Chris Heyde in his office at Columbia University with interviewers Steven Kou (left) and Paul Glasserman (right) (December 2005).*

**Paul:** What are your research plans for the future?

**Chris:** I am one of the Chief Investigators of the newly established Australian Research Council "Centre of Excellence for Mathematics and Statistics of Complex Systems" and most of the research that I am doing and planning is associated with the themes of this center. This includes the nonstandard limit theorems that occur when long-range dependence holds and the classical theory breaks down. I have had a long-term interest in linking the physical explanations with the mathematical ones for such phenomena as long-range dependence and intermittency, and for exploring fractal behavior and scaling properties such as given by self-similarity. These topics are quite closely associated with diverse applications, for example, in risky asset and teletraffic modeling, and in the earth and environmental sciences. There is plenty to keep me busy.